\documentclass{goose-article}

\usepackage{algorithm,algpseudocode,setspace}
\usepackage[procnames]{listings}
\usepackage{color}
\usepackage{soul}
\usepackage{verbatim}
\usepackage{xcolor}

\usepackage[normalem]{ulem}

\definecolor{keywords}{RGB}{255,0,90}
\definecolor{comments}{RGB}{0,0,113}
\definecolor{red}{RGB}{160,0,0}
\definecolor{green}{RGB}{0,150,0}

\catcode`@=11
\def\nvphantom{\v@true\h@false\nph@nt}
\def\nhphantom{\v@false\h@true\nph@nt}
\def\nphantom{\v@true\h@true\nph@nt}
\def\nph@nt{\ifmmode\def\next{\mathpalette\nmathph@nt}%
  \else\let\next\nmakeph@nt\fi\next}
\def\nmakeph@nt#1{\setbox\z@\hbox{#1}\nfinph@nt}
\def\nmathph@nt#1#2{\setbox\z@\hbox{$\m@th#1{#2}$}\nfinph@nt}
\def\nfinph@nt{\setbox\tw@\null
  \ifv@ \ht\tw@\ht\z@ \dp\tw@\dp\z@\fi
  \ifh@ \wd\tw@-\wd\z@\fi \box\tw@}

\newcommand{\T}[1]{\bm{#1}}
\newcommand{\TT}[1]{\mathbb{#1}}
\newcommand{\mat}[1]{\underline{\underline{#1}}\vphantom{#1}}
\newcommand{\matT}[1]{\underline{\underline{\bm{#1}}}\vphantom{\bm{#1}}}
\newcommand{\matTT}[1]{\underline{\underline{\mathbb{#1}}}\vphantom{\mathbb{#1}}}
\newcommand{\col}[1]{\underline{#1}\vphantom{#1}}
\newcommand{\colT}[1]{\underline{\bm{#1}}\vphantom{\bm{#1}}}

\hypersetup{pdfauthor={T.W.J. de Geus}}

\title{%
  Finite strain FFT-based non-linear solvers made simple
}
\author[1,2,5]{T.W.J.~de~Geus}
\author[3]{J.~Vond\v{r}ejc}
\author[4]{J.~Zeman}
\author[2]{R.H.J.~Peerlings}
\author[2]{M.G.D.~Geers}
\affil[1]{
  Materials Innovation Institute (M2i),\nl
  P.O.~Box~5008, 2600~GA~Delft, The~Netherlands
}
\affil[2]{
  Eindhoven University of Technology, Department of Mechanical Engineering,\nl
  P.O.~Box~513, 5600~MB~Eindhoven, The~Netherlands
}
\affil[3]{
  Technische Universit\"{a}t Braunschweig, Institute of Scientific Computing,\nl
  D-38092 Braunschweig, Germany
}
\affil[4]{
  Czech Technical University in Prague, Department of Mechanics, Faculty of Civil Engineering,\nl
  Th\'{a}kurova 7, 166 29 Prague 6, Czech Republic
}

\contact{%
  $^*$Corresponding author now at \'{E}cole polytechnique f\'{e}d\'{e}rale de Lausanne (EPFL), Station 18, 1015 Lausanne, Switzerland: \href{mailto:tom.degeus@epfl.ch}{tom.degeus@epfl.ch}, \href{mailto:tom@geus.me}{tom@geus.me} %
}

\header{%
  T.W.J.~de~Geus, J.~Vond\v{r}ejc, J.~Zeman, Peerlings, R.H.J.~Peerlings, M.G.D.~Geers\\ Comput.\ Methods Appl.\ Mech.\ Eng., 2017, 318:412–430, \doi{10.1016/j.cma.2016.12.032}, \eprint{1603.08893}
}

\begin{document}

\maketitle

\begin{abstract}
Computational micromechanics and homogenization require the solution of the mechanical equilibrium of a periodic cell that comprises a (generally complex) microstructure. Techniques that apply the Fast Fourier Transform have attracted much attention as they outperform other methods in terms of speed and memory footprint. Moreover, the Fast Fourier Transform is a natural companion of pixel-based digital images which often serve as input. In its original form, one of the biggest challenges for the method is the treatment of (geometrically) non-linear problems, partially due to the need for a uniform linear reference problem. In a geometrically linear setting, the problem has recently been treated in a variational form resulting in an unconditionally stable scheme that combines Newton iterations with an iterative linear solver, and therefore exhibits robust and quadratic convergence behavior. Through this approach, well-known key ingredients were recovered in terms of discretization, numerical quadrature, consistent linearization of the material model, and the iterative solution of the resulting linear system. As a result, the extension to finite strains, using arbitrary constitutive models, is at hand. Because of the application of the Fast Fourier Transform, the implementation is substantially easier than that of other (Finite Element) methods. Both claims are demonstrated in this paper and substantiated with a simple code in Python of just 59 lines (without comments). The aim is to render the method transparent and accessible, whereby researchers that are new to this method should be able to implement it efficiently. The potential of this method is demonstrated using two examples, each with a different material model.
\end{abstract}


\section{Introduction}

Computational micromechanics and homogenization generally involve the numerical solution of the mechanical equilibrium of a periodic unit-cell \citep{Michel1999,Kanoute2009,Geers2010,Bensoussan1978}. Such a unit-cell thereby provides a representative geometrical representation of the microstructure -- which is often complex. An accurate representation of reality therefore necessitates a high-resolution numerical method, which remains efficient in three dimensions. The conventional approach, using Finite Elements, results in costly computations. For periodic cells, an attractive competitor to the Finite Element Method was developed by \citet{Moulinec1994,Moulinec1998}. It employs the Fast Fourier Transform (FFT) to obtain a significant gain in efficiency compared to Finite Elements, both in terms of speed and in terms of memory footprint. Furthermore, the method is rather straightforward to implement efficiently, as FFT-libraries are readily available. In its original form, the method uses a homogeneous auxiliary elastic problem. This form has been extended to (geometrically) non-linear problems \citep{Michel2001,Vinogradov2008,Monchiet2012,Kabel2014,Eisenlohr2013,Shanthraj2015}, however still relying on a reference medium, which may impact the convergence and obscures consistent linearization. It has nevertheless been successfully applied to large-scale computations \citep{Diehl2015,Robert2015,Montagnat2014,Lebensohn2008}. Currently, even a mature open-source implementation is available featuring different solution techniques \citep{damask}. It does however not include consistent linearization, which is the contribution of the present paper.

Recently, \citet{Vondrejc2013} have recognized that the problem can be reformulated in a variational form. This framework has the benefit that discretization, quadrature, constitutive linearization, and the solution of a linear system can be properly distinguished and optimized individually \citep{Moulinec2013,Mishra2015}. This shows that the equations solved are essentially nodal equilibrium equations -- much like in the Finite Element Method (FEM). The variational formulation also obsoletes an auxiliary elastic problem by employing a projection operator that merely maps arbitrary tensor fields to compatible ones. Unlike conventional FEM, the integration points coincide with the nodes, which results in an extremely efficient, local, update of the strains and stresses. The method allows one to use the local consistent tangent in the equilibrium iterations, providing a quadratic convergence rate. It has recently been demonstrated that convergence is robust for several types of non-linear constitutive behavior in small strains \citep{Zeman2016}.

Treating the problem in the variational form renders the extension to finite strains relatively straightforward. The intrinsic simplicity of the method is thereby reconfirmed. This claim is demonstrated in the present (short) paper, by incorporating a simple Python code for a three-dimensional unit-cell problem in finite strains, which enables efficient and independent implementation by other researchers. The computational efficiency is inherited from the efficient Fast Fourier Transform. We furthermore demonstrate that the method can be employed for arbitrary constitutive models in finite strains, defined in either the undeformed or the deformed configuration. For the latter, both the stress and the consistent tangent must be transformed to the undeformed configuration using a pull-back involving simple tensor operations. Again the simplicity of the method results in a compact Python code of well below 200 lines. All codes have been made available for free downloading \citep{GooseFFT} whereby we would like to invite the interested reader to contribute with other relevant examples.

This paper is structured as follows. In Section~\ref{sec:method} the variational approach due to \citet{Vondrejc2013} and \citet{Zeman2016} is extended to finite strains, following the same steps. Section~\ref{sec:implement} presents a pseudo-algorithm accompanied with a discussion on the most import operations. In Sections~\ref{sec:elas} and~\ref{sec:simo}, two examples are considered. A short summary is included in Section~\ref{sec:summary}. The adopted nomenclature is presented in Appendix~\ref{sec:nomenclature}, the Python code in Appendix~\ref{sec:code:python}, and the treatment of even-sized grids in Appendix~\ref{sec:non-odd-grid} (odd-sized grids are used in the main text). Finally the projection operator, which distinguishes this method from FEM and other FFT-based methods, is explained in more detail in Appendix~\ref{sec:projection}.

\section{Method}
\label{sec:method}

\subsection{Problem statement}

\begin{figure}[htp]
  \centering
  \includegraphics[width=.3\textwidth]{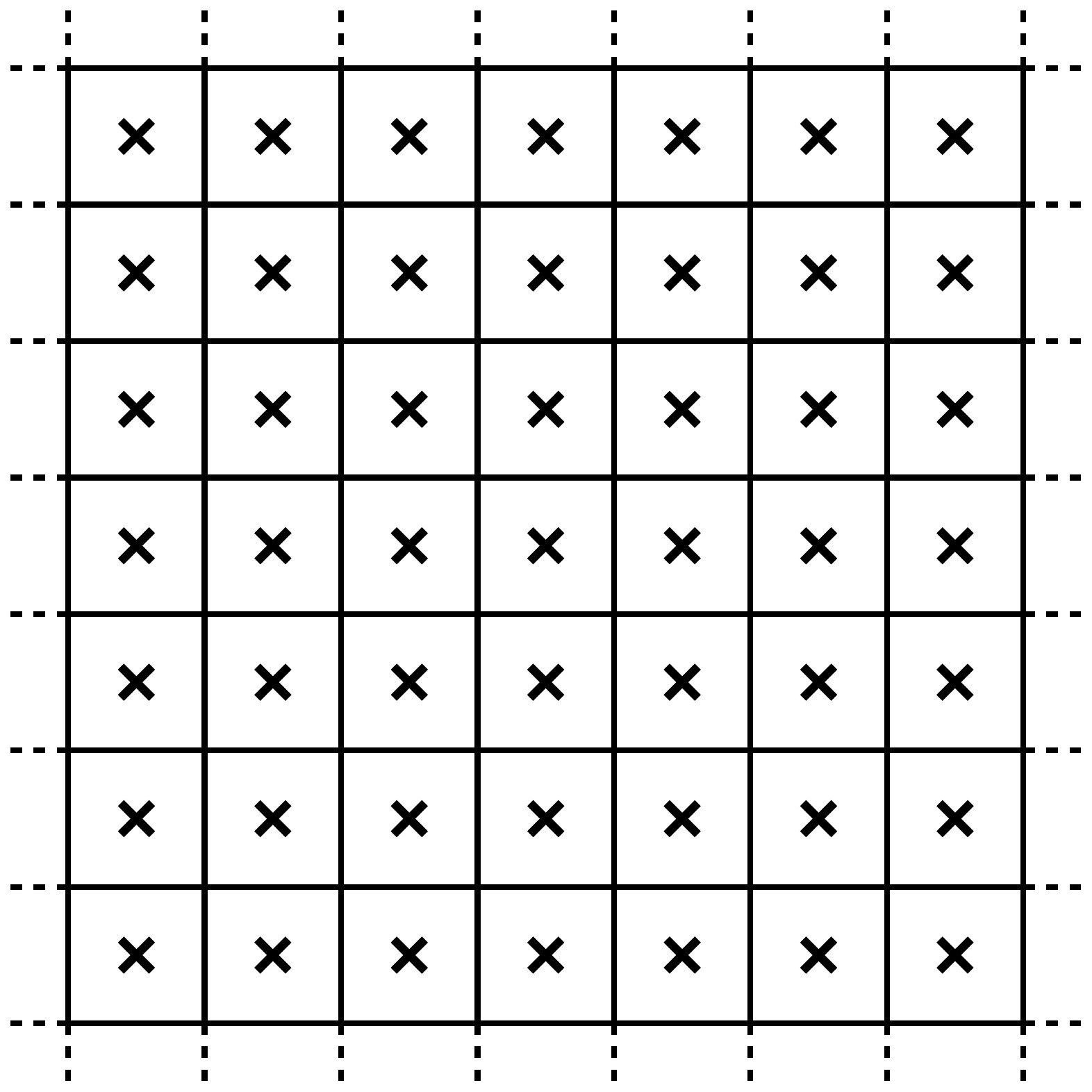}
  \caption{A simple periodic cell discretized using a regular pixel grid, each `pixel' having equal area. The solution is discretized on the same grid (Eq.~\eqref{eq:shape:F}); the corresponding nodes are shown using crosses. The nodes and integration points coincide. }
  \label{fig:method:cell}
\end{figure}

A periodic cell is considered, containing a representative volume of the microstructure (in one, two, or three dimensions). It is composed of one or more phases described by arbitrary material models formulated in a finite strain framework. For simplicity of notation, the derivations in this section are performed for a non-linear elastic model. The extension to time- or history-dependent materials is straightforward (in line with FEM-based formulations~\citep{Zeman2016}) and therefore only included as pseudo-algorithm and in the accompanying examples \citep{GooseFFT}. The microstructure is discretized using a regular grid (i.e.\ a grid of pixels or voxels, see Figure~\ref{fig:method:cell}).

The goal is to solve for static mechanical equilibrium in the periodic cell for a given applied overall deformation. The balance of linear momentum, pulled-back to the (undeformed) reference configuration, reads
\begin{equation}\label{eq:divP}
  \vec{\nabla}_0 \cdot \T{P}^T = \vec{0}
\end{equation}
involving the divergence with respect to the reference configuration of the transposed first Piola-Kirchhoff stress tensor $\T{P}$. In index notation Eq.~\eqref{eq:divP} reads $\partial P_{ij} / \partial X_j$, see Appendix~\ref{sec:nomenclature} for our nomenclature. The stress $\T{P}$ depends non-linearly on the deformation gradient $\T{F}$:
\begin{equation}\label{eq:P-F}
  \T{P} = \T{P} ( \T{F} )
\end{equation}
%

\subsection{Weak form}

The integral form is obtained by multiplying \eqref{eq:divP} with test functions $\delta \vec{x}$ and integrating over the reference domain $\Omega_0$:
\begin{equation}
  \int\limits_{\Omega_0}
    \delta \vec{x} \cdot \big( \vec{\nabla}_0 \cdot \T{P}^T \big)
  \; \mathrm{d} \Omega_0
  =
  0
\end{equation}
which must hold for all periodic $\delta \vec{x}$.

Subsequently, integration by parts is applied in conjunction with Gauss' divergence theorem. The boundary term that arises vanishes because of periodicity. The result reads
\begin{equation}\label{eq:weak:x}
  \int\limits_{\Omega_0}
    \big( \vec{\nabla}_0 \, \delta \vec{x} \big)^T : \T{P}^T
  \; \mathrm{d} \Omega_0
  =
  0
\end{equation}
where $:$ denotes a double tensor contraction, see again Appendix~\ref{sec:nomenclature}.

In Finite Elements, the problem is discretized at this point by introducing a finite element interpolation scheme for the position vector $\vec{x}(\vec{X})$ -- or the displacement $\vec{u}(\vec{X})$ -- as well as the virtual displacement $\delta \vec{x}(\vec{X})$, where $\vec{X}$ is the position in the undeformed configuration. The problem is then solved for the nodal positions or displacements. In the FFT-methods, one generally does not formulate the problem in terms of positions/displacements but in terms of strains, or, in the context of finite deformation, in terms of the deformation gradient tensor $\T{F}$. By analogy, the test function is expressed in terms of a virtual deformation gradient. The weak form of Eq.~\eqref{eq:weak:x} is accordingly reformulated as
\begin{equation}\label{eq:weak:F:incomp}
  \int\limits_{\Omega_0}
    \delta \T{F} : \T{P}^T
  \; \mathrm{d} \Omega_0
  =
  0
\end{equation}
in which the test functions $\delta \T{F}$ are periodic and compatible. Note that compatibility is guaranteed when $\delta \T{F}$ is the gradient of a virtual position vector, as in Finite Elements, but now must be enforced as a constraint in conjunction with Eq.~\eqref{eq:weak:F:incomp}.

\subsection{Projection to a compatible solution space}

The compatibility of the test functions $\delta \T{F}$ is imposed by means of a projection operator $\TT{G}$. It maps an arbitrary field $\tilde{\T{A}}$ to its compatible part $\T{A}$ through
\begin{equation}\label{eq:compatibility}
  \TT{G} \star \tilde{\T{A}} = \T{A}
\end{equation}
wherein $\star$ is the convolution operator. The convolution can be evaluated in Fourier space as a simple, local, double tensor contraction. Furthermore, $\TT{G}$ has a simple closed-form expression in Fourier space, see below (Eq.~\eqref{eq:G}). Its background and interpretation is discussed in Appendix~\ref{sec:projection}.

Application of Eq.~\eqref{eq:compatibility} to the weak form of Eq.~\eqref{eq:weak:F:incomp} results in
\begin{equation}\label{eq:weak:F}
  \int\limits_{\Omega_0}
    ( \TT{G} \star \delta \tilde{\T{F}} ) : \T{P}^T
  \; \mathrm{d} \Omega_0
  =
  \int\limits_{\Omega_0}
    \delta \tilde{\T{F}} : ( \TT{G} \star \T{P} )
  \; \mathrm{d} \Omega_0
  =
  0
\end{equation}
whereby the symmetry of $\TT{G}$ has been used. Equation~\eqref{eq:weak:F} should now hold for arbitrary, i.e.\ not necessarily compatible, periodic test functions $\delta \tilde{\T{F}}$. Please note that the deformation gradient $\T{F}$, hidden in the stress $\T{P}$ through Eq.~\eqref{eq:P-F}, should still satisfy the compatibility constraint. This is enforced below, in Section~\ref{sec:method:lin}.

\subsection{Discretization}

Adopting a Galerkin scheme, the unknown field $\T{F}$ and the test functions $\delta \tilde{\T{F}}$ are discretized in the same way. Like in Finite Elements, the continuous fields $\T{F}$ and $\delta \tilde{\T{F}}$ are approximated by a finite number of $n$ nodal values, $\T{F}_{k}$ and $\delta \tilde{\T{F}}_{k}$, that are multiplied with shape functions $N_{k}$ associated with each node, i.e.
\begin{alignat}{2}
  \T{F} (\vec{X})
  &\approx
  \sum\limits_{k=1}^n N_{k} ( \vec{X} ) \, \phantom{\delta} \T{F}_{k} \;
  &=
  \col{N}^\mathsf{T} (\vec{X}) \; \phantom{\delta} \colT{F}
  \label{eq:shape:F}
  \\
  \delta \tilde{\T{F}} (\vec{X})
  &\approx
  \sum\limits_{k=1}^n N_{k} ( \vec{X} ) \, \delta \tilde{\T{F}}_{k} \;
  &=
  \col{N}^\mathsf{T} (\vec{X}) \; \delta \tilde{\colT{F}}
  \label{eq:shape:deltaF}
\end{alignat}
where the underlined symbols indicate a column matrix of nodal quantities and the superscript $\mathsf{T}$ transposition of the matrix, e.g.:
\begin{equation}
  \colT{F} =
  \Big[
    \T{F}_{1} ,\;\;
    \T{F}_{2} ,\;\;
    \ldots    ,\;\;
    \T{F}_{n}
  \Big]^\mathsf{T}
\end{equation}
In Figure~\ref{fig:method:cell} the nodes are indicated by crosses. Fundamental trigonometric polynomials are used as shape functions $\col{N}$ (see \citep{Saranen2002}). These functions have global support, but are interpolatory, i.e.\
\begin{equation}\label{eq:shape:delta}
  N_{k} (\vec{X}_m) = \delta_{km}
\end{equation}
where $\delta_{km}$ is the Kronecker delta and $k$ and $m$ are node numbers. In addition, they constitute a partition of unity, i.e.\ the sum of all shape functions equals one everywhere in the domain.

The discretization is applied to the weak form (Eq.~\eqref{eq:weak:F}) which therefore becomes
\begin{equation}\label{eq:weak:F:discr}
  \Big( \delta \tilde{\colT{F}} \Big)^\mathsf{T} :
  \int\limits_{\Omega_0}
    \col{N} \; ( \TT{G} \star \T{P} )
  \; \mathrm{d} \Omega_0
  =
  0
\end{equation}
where the fact that the nodal quantities $\delta \tilde{\T{F}}_k$ are independent of $\vec{X}$ has been used to take them out of the integral. Note that in this discretized form, the stress $\T{P}$ depends on $\col{N}^\mathsf{T} \colT{F}$.

Eq.~\eqref{eq:weak:F:discr} must hold for all $\delta \tilde{\colT{F}}$, which implies
\begin{equation}\label{eq:weak:discr}
  \int\limits_{\Omega_0}
    \col{N} \; ( \TT{G} \star \T{P} )
  \; \mathrm{d} \Omega_0
  =
  \colT{O}
\end{equation}
%

\subsection{Quadrature}

The integration is performed using numerical quadrature. For the employed trigonometric polynomials $\col{N}$ the trapezoidal rule is used, which assigns equal weight to each nodal quantity. This implies that there is no distinction between integration points and nodes. Applied to Eq.~\eqref{eq:weak:discr},
\begin{equation}\label{eq:weak:quad}
  \sum\limits_k \;
  \col{N} ( \vec{X}_k ) \; \big[ \TT{G} \star \T{P} \big] (\vec{X}_k)
  =
  \colT{O}
\end{equation}
We now exploit the fact that the shape functions can be expressed in terms of the discrete Fourier coefficients; and the delta property of the shape functions (Eq.~\eqref{eq:shape:delta}). In terms of the Discrete Fourier Transform $\mathcal{F}$ and its inverse $\mathcal{F}^{-1}$, Eq.~\eqref{eq:weak:quad} can be expressed as
\begin{equation}\label{eq:weak:discr:FFT-intermediate}
  \mathcal{F}^{-1}
  \big\{\; \matTT{\hat{G}} : \mathcal{F} \{ \colT{P} \}\; \big\}
  =
  \colT{O}
\end{equation}
see Ref.~\citep{Zeman2016} for details. It is thus composed of a sequence of (i) the Fourier transform $\mathcal{F}$ of the nodal stresses $\colT{P}$, (ii) a double tensor contraction of them with the, explicitly known and constant, Fourier-coefficients $\matTT{\hat{G}}$, and (iii) the inverse Fourier transform $\mathcal{F}^{-1}$ of the final result. The notation of Eq.~\eqref{eq:weak:quad} is abbreviated to
\begin{equation}\label{eq:weak:discr:FFT}
  \matTT{G} : \colT{P}
  =
  \colT{O}
\end{equation}
Thereby, the stress at each node depends locally on the deformation gradient, i.e.
\begin{equation}
  \colT{P} =
  \Big[
    \T{P}_1 ( \T{F}_1 ) ,\;\;
    \T{P}_2 ( \T{F}_2 ) ,\;\;
    \ldots              ,\;\;
    \T{P}_n ( \T{F}_n )
  \Big]^\mathsf{T}
\end{equation}
Equation~\eqref{eq:weak:discr:FFT} can thus be thought of as a set of local, tensor valued, constitutive relations of the form $\T{P}_k (\T{F}_k)$, whereby the global equilibrium is enforced by the non-locality of the projection operator $\matTT{G}$.

In Fourier space, the matrix $\matTT{\hat{G}}$ is a diagonal matrix,\footnote{N.B. the matrix is diagonal in terms of tensors, it is therefore block diagonal in terms of tensor components.} i.e.
\begin{equation}
  \matTT{\hat{G}} =
  \begin{bmatrix}
    \hat{\TT{G}}_{11} &                   &        &                   \\
                      & \hat{\TT{G}}_{22} &        &                   \\
                      &                   & \ddots &                   \\
                      &                   &        & \hat{\TT{G}}_{nn}
  \end{bmatrix}
\end{equation}

The expression for $\hat{\TT{G}}_{kk}$, restricted to grids with an odd number of nodes, reads
\begin{equation}\label{eq:G}
(\hat{G}_{kk})_{ijlm} (\vec{q}_k) =
\begin{cases}
  \displaystyle
  0
  & \displaystyle\vphantom{\frac{\vec{\xi}}{\vec{\xi}}}
  \text{for }\vec{q}_k = \vec{0}
  \\
  \displaystyle
  \frac{
    \delta_{im}\, \xi_j(\vec{q}_k)\, \xi_l(\vec{q}_k)
  }{
    \| \vec{\xi} \|^2
  }
  & \text{otherwise}
\end{cases}
\end{equation}
wherein $\vec{q}_k$ is the (spatial) frequency vector and $\xi_i$ are the scaled frequencies that account for the size of the cell through $\xi_i (\vec{q}) = q_i / L_i$ (with $L_i$ the size of the periodic cell in direction $i$). The zero frequency ($\vec{q} = \vec{0}$) is associated with the mean of the inversely transformed field. The definition in Eq.~(\ref{eq:G}a) therefore ensures the zero mean property
\begin{equation}\label{eq:G:mean}
  \int\limits_{\Omega_0} \matTT{G} : \colT{P} \; \mathrm{d}\Omega_0 = \T{O}
\end{equation}
This is used below to prescribe the macroscopic deformation, while solving only for the periodic micro-fluctuations. For grids with an even number of nodes, enforcing the zero-mean property of \eqref{eq:G:mean} along with compatibility is slightly more involved, see Appendix~\ref{sec:non-odd-grid}.

\subsection{Linearization}
\label{sec:method:lin}

The weak form in Eq.~\eqref{eq:weak:discr:FFT} is a non-linear equation, as the material model involves a non-linear relation between the first Piola-Kirchhoff stress and the deformation gradient. We employ Newton iterations to solve the nodal equilibrium equations~\eqref{eq:weak:discr:FFT}. To this end the nodal unknowns at iteration $i+1$ are expressed as
\begin{equation}\label{eq:iter:F}
  \colT{F}_{(i+1)} = \colT{F}_{(i)} + \delta \colT{F}
\end{equation}
where $\colT{F}_{(i)}$ are the last known iterative values of the deformation gradients and $\delta \colT{F}$ are their iterative updates. Note that $\delta$ is now used to indicate a small variation. The stresses are linearized around $\colT{F}_{(i)}$. In a material point this corresponds to
\begin{equation}\label{eq:iter:varP}
  \delta \T{P}^T = \TT{K}_{(i)} : \delta \T{F}^T
\end{equation}
Note that by expressing the linearization in this way, including the transposes, we stay as close as possible to the standard Finite Element formulation.

Combined with the discretized weak form in Eq.~\eqref{eq:weak:discr:FFT}, the iterative update $\delta \colT{F}$ is found by solving the following linearized system
\begin{equation}\label{eq:iter:system}
  \matTT{G} : \matTT{K}_{(i)}^{LT} : \delta \colT{F}^T
  = -
  \matTT{G} : \colT{P}_{(i)}
\end{equation}
wherein $\matTT{K}_{(i)}^{LT}$ are the left-transposed local tangent stiffnesses, assembled in a diagonal matrix:\footnote{N.B. the matrix is diagonal in terms of tensors, it is therefore block diagonal in terms of tensor components.}
\begin{equation}\label{eq:mat:Kmat}
  \matTT{K}_{(i)}^{LT} =
  \begin{bmatrix}
    \big(\TT{K}_{(i)}^{LT}\big)_{11} & & \\
    & \big(\TT{K}_{(i)}^{LT}\big)_{22} & \\
    & & \ddots & \\
    & & & \big(\TT{K}_{(i)}^{LT}\big)_{nn}
  \end{bmatrix}
\end{equation}
It thus consists only of the constitutive tangents evaluated locally at the nodes (Eq.~\eqref{eq:iter:varP}).

Note that the compatibility of the deformation gradient field still needs to be enforced. This is done by solving the linear system of Eq.~\eqref{eq:iter:system} using an iterative solver which delivers a compatible solution in each iteration \cite{Vondrejc2013,Vondrejc2016,Zeman2016,Moulinec2013,Mishra2015}. To satisfy compatibility during the entire iterative process, projection based iterative methods such as e.g.\ the conjugate gradient (CG) method and the generalized minimal residual method (GMRES), Chebyshev iterations, or Richardson iterations (used in the original Moulinec-Suquet algorithm \cite{Moulinec1994,Moulinec1998}) can be used \cite{Mishra2015,Moulinec2013}. Alternatively, compatibility is satisfied only at convergence using the accelerated method of Ref.~\cite{Eyre1999}.

\subsection{Connection to existing FFT-based solvers}
\label{sec:connections}

The main difference between the state-of-the art finite-strain solvers~\cite{Kabel2014,Shanthraj2015} and the present formulation lies in the application and interpretation of the projection operator. In our approach, the projection operator $\matTT{G}$ merely ensures that each row of the compatible field $\T{A}$ in Eq.~\eqref{eq:compatibility} can be obtained as the gradient of a scalar potential; see Appendix~\ref{sec:projection} for further clarification. This construction is obviously independent of the constitutive laws used. It therefore requires no parameters or other problem specific choices, and $\matTT{G}$ remains constant throughout the entire simulation.

In contrast, the formulations introduced in Refs.~\cite{Kabel2014,Shanthraj2015} enforce the compatibility by means of a Green operator,  associated with a reference problem with an auxiliary constant stiffness $\TT{C}_0$. However, the choice of the reference stiffness in the finite-strain setting is not obvious and $\TT{C}_0$ is mostly introduced on a heuristic basis. For example, \citet{Kabel2014} suggest to determine $\TT{C}_0$ from the extreme positive eigenvalues of local tangent stiffnesses (i.e.\ from $(\TT{K}_{(i)})_{11}$ through $(\TT{K}_{(i)})_{nn}$ in Eq.~\eqref{eq:mat:Kmat}), whereas \citet{Shanthraj2015} employ a weighted average of the local tangents. In both cases, the Green operator must be updated at every Newton iteration, which can become costly.

\section{Implementation}
\label{sec:implement}

The numerical algorithm requires the solution of \eqref{eq:weak:discr:FFT} in an incremental--iterative fashion. Each increment thereby consists of Newton iterations updating the nodal deformation gradients $\colT{F}_{(i+1)}$ using (\ref{eq:iter:F}--\ref{eq:iter:system}) until equilibrium is satisfied up to an accuracy $\eta^\mathrm{NW}$, by employing the linearized constitutive response. The linear system in Eq.~\eqref{eq:iter:system} is solved up to an accuracy $\eta^\mathrm{CG}$ using the conjugate gradient iterative solver.

\subsection{Boundary conditions}

With the periodic micro-fluctuations of $\colT{F}$ following from equilibrium and compatibility, only the macroscopic deformation or stress needs to be prescribed. In this paper we restrict ourselves to a fully prescribed macroscopic deformation gradient $\bar{\T{F}}$, as this is the easiest and the most efficient choice. More general cases are discussed in the literature \cite[e.g.][]{Kabel2015b}.

We start from an equilibrium state given by $\colT{F}_{(0)}$, and apply a macroscopic deformation gradient $\bar{\T{F}}$. More specifically, we apply the difference of $\bar{\T{F}}$ and to the mean of $\colT{F}_{(0)}$:
\begin{equation}
  \Delta \bar{\T{F}}
  =
  \bar{\T{F}} - \int_{\Omega_0} \colT{F}_{(0)} \; \mathrm{d}\Omega_0
\end{equation}

For the first Newton iteration, equilibrium reads
\begin{equation} \label{eq:bc:nonlin}
  \matTT{G} :
  \colT{P} \big( \colT{F}_{(0)} + \Delta \bar{\colT{F}} + \delta \colT{F} \big)
  =
  \colT{O}
\end{equation}
where $\Delta \bar{\colT{F}}$ is a column with $\Delta \bar{\T{F}}$ on each row. Linearization of \eqref{eq:bc:nonlin} results in
\begin{equation} \label{eq:bc:lin}
  \matTT{G} : \matTT{K}_{(0)}^{LT} : \delta \colT{F}^T
  = -
  \matTT{G} : \matTT{K}_{(0)}^{LT} : \Delta \bar{\colT{F}}^T
\end{equation}
where $\matTT{K}_{(0)}$ is the constitutive tangent about $\colT{F}_{(0)}$. Note that use has been made of the fact that $\colT{F}_{(0)}$ is in equilibrium, i.e.\ that $\matTT{G} : \colT{P} \big( \colT{F}_{(0)} \big) = \colT{O}$. After solving the system in Eq.~\eqref{eq:bc:lin} we set
\begin{equation}
  \colT{F}_{(1)}
  =
  \colT{F}_{(0)} + \Delta \bar{\colT{F}} + \delta \colT{F}
\end{equation}
and proceed as normally (Eqs.~(\ref{eq:iter:F}--\ref{eq:iter:system})). It is thereby important to point out that the definition of $\matTT{G}$ (Eq.~\eqref{eq:G}) ensures that the mean
\begin{equation}\label{eq:mean:zero}
  \int_{\Omega_0} \delta \colT{F} \; \mathrm{d}\Omega_0
  =
  \T{O}
\end{equation}
All iterations will thus satisfy the prescribed $\bar{\T{F}}$ exactly.

The interpretation of Eq.~\eqref{eq:bc:lin} is that, by solving the linear system, the macroscopic deformation, $\Delta \bar{\T{F}}$, is distributed over this microstructure using the tangent $\matTT{K}_{(0)}$, which contains the microstructural heterogeneity. Eq.~\eqref{eq:bc:lin} therefore has strong similarities with the application of essential (Dirichlet) boundary conditions in the Finite Element Method.

\subsection{Program structure}

The above implementation aspects are further explained using Algorithm~\ref{fig:algorithm}, which holds for the case that the macroscopic deformation gradient tensor $\bar{\T{F}}$ is prescribed. As the goal of this paper is to make the method \textit{simple}, an implementation, in the Python programming language, is included in Appendix~\ref{sec:code:python} (for the hyper-elastic example given in Section~\ref{sec:elas}). It depends only on standard Python and its scientific libraries, i.e.\ no custom software or libraries are used. Note that this code, and all other codes used for this paper are freely available for download and use \citep{GooseFFT}, and that we invite similar contributions there.
\begin{algorithm}[htp]
\newcommand{\rcomment}[1]{\Comment\normalsize#1\normalsize}
\caption{}
\label{fig:algorithm}
\setstretch{1.3}
\begin{algorithmic}[1]
\vspace*{1em}

  \State
  Initialize $\colT{F}_{(0)} = \colT{I}$ and history variables

  \vspace*{0.9em}
  \For{$ t = 0 , \Delta t , 2 \Delta t , ... $} \rcomment{\textbf{\underline{incremental loop}}}

    \vspace*{0.9em}
    \For{$ i = 0, 1, 2, ... $}
    \rcomment{\textbf{\underline{Newton iterations}}}

      \vspace*{0.9em}
      \State
      \(
      \colT{F}_{(i)} \rightarrow \matTT{K}_{(i)}, \colT{P}_{(i)}
      \)
      \rcomment{\emph{constitutive response}}

      \vspace*{0.9em}
      \If{\( i=0 \)} \rcomment{\emph{boundary condition}}

        \State
        Solve
        \(
        \matTT{G} : \matTT{K}_{(i)}^{LT} : \delta \colT{F}^T = -
        \matTT{G} : \matTT{K}_{(i)}^{LT} : \Delta \bar{\colT{F}}^T
        \)
        with accuracy $\eta^\mathrm{CG}$
        \State
        \(
        \colT{F}_{(i+1)} =
        \colT{F}_{(i)}   + \delta \colT{F} + \Delta \bar{\colT{F}}
        \)

      \Else \rcomment{\emph{equilibrium iteration}}

        \State
        Solve
        \(
        \matTT{G} : \matTT{K}_{(i)}^{LT} : \delta \colT{F}^T = -
        \matTT{G} : \colT{P}_{(i)}
        \)
        with accuracy $\eta^\mathrm{CG}$
        \State
        \(
        \colT{F}_{(i+1)} =
        \colT{F}_{(i)}   + \delta \colT{F}
        \)

      \EndIf

      \vspace*{0.9em}
      \If{\(
        \|\, \delta \colT{F} \,\| \,/\,
        \|\, \colT{F}_{(t)} \,\| > \eta^\text{NW}
        \;\mathrm{and}\; i > 0
      \) }
      \rcomment{\emph{convergence criterion}}

      \State
      Proceed to l.\ 16 (end of iterative loop)

      \EndIf

    \vspace*{0.9em}
    \EndFor

    \vspace*{0.9em}
    \State
    \( \colT{F}_{(t + \Delta t)} = \colT{F}_{(i+1)} \)
    \rcomment{\emph{store converged state}}
    \State
    Store history variables

  \vspace*{0.9em}
  \EndFor

\vspace*{1em}
\end{algorithmic}
\setstretch{1.0}
\end{algorithm}

\subsection{Implementation aspects}

Based on Algorithm~\ref{fig:algorithm} the following comments are made:
\begin{itemize}
  \item The conjugate gradient solver only requires the (repeated) result of $\matTT{G} : \matTT{K}_{(i)}^{LT} : \delta \colT{F}^T$. Therefore $\matTT{\hat{G}}$ and $\matTT{K}_{(i)}$ are stored independently, the fully populated system matrix is never assembled.
  \item Although $\colT{F}$, $\colT{P}$, $\matTT{K}$, and $\matTT{\hat{G}}$ are written as columns and diagonal matrices, the actual storage coincides with the grid of nodes. For a three-dimensional problem, the grids of second-order tensors are stored as five-dimensional matrices with (only) $3^2 \times n$ entries, where $3^2$ corresponds to the tensor components and $n$ to the nodes. Likewise, the fourth-order tensors are stored as seven-dimensional matrices with $3^4 \times n$ entries. More involved representations (e.g.\ sparse matrices) are therefore not needed.
  \item The projection $\matTT{G} : \bullet$ involves the Fourier transform, a tensor contraction, and the inverse Fourier transform. Both the Fourier transform, $\mathcal{F}$, and the inverse Fourier transform, $\mathcal{F}^{-1}$, of a second-order tensor are performed per tensor component for each spatial direction. For a three-dimensional problem it thus involves $3 \times 3 \times 3$ one-dimensional Fourier transforms and the same number of inverse Fourier transforms. Most FFT-libraries however provide an interface to such an operation in a single command.
  \item For most problems, the Fourier transform and its inverse are the most costly operations. Fortunately, efficient open-source implementations with interfaces to all popular programming languages are freely available, often readily parallelized. Furthermore, the Fourier transform only scales as $n \log (n)$ with the problem size $n$. A practical disadvantage can be that the efficiency is highly dependent on the exact number of grid points. A speed-up of orders of magnitude can sometimes be obtained by slightly modifying the number of grid points, whereby the optima are found when $n$ is a power of two \citep{Frigo2005}. Not all problems allow for this luxury.
  \item The constitutive model is evaluated only at the nodal level, requiring only the nodal deformation gradient tensor $\T{F}_k$ as input. The constitutive implementation at the integration point level as used in Finite Element codes can therefore directly be used here as well.
\end{itemize}
%

\section{Application to Hyper-elasticity, formulated in the reference configuration}
\label{sec:elas}

\subsection{Introduction}

One of the simplest constitutive models to consider is a hyper-elastic model defined in the reference configuration. It is therefore used for the implementation in Appendix~\ref{sec:code:python}, and was also used by \citet{Kabel2014}.\footnote{Note that this model is not a proper elastic constitutive model for very large deformations (see also \citep{Bertram2007}) It nevertheless serves the purpose of being simple and sufficiently realistic within the range of deformations analyzed here. Furthermore it allows direct comparison with the analytical solution to a 2-D laminate problem by \citet{Kabel2014}, for which a perfect agreement has been found using our presented implementation (comparison not shown).}. Note that the problem is still non-linear because of the geometric non-linearity.

\subsection{Constitutive model}

The model is defined in the undeformed configuration and it involves a linear relation between the second Piola-Kirchhoff stress $\T{S}$ and the Green-Lagrange strain $\T{E}$:
\begin{equation}\label{eq:elas:S}
  \T{S} = \TT{C} : \T{E}
\end{equation}
wherein $\TT{C}$ is the standard fourth-order elastic stiffness tensor
\begin{equation}\label{eq:elas:C}
  \TT{C}
  =
  \lambda \, \T{I} \otimes \T{I}
  +
  2 \mu \; \TT{I}^s
\end{equation}
with Lam\'{e}'s constants $\lambda$ and $\mu$. In terms of Young's modulus $E$ and Poisson's ratio $\nu$ they read:
\begin{equation}
  \lambda = \frac{E \nu}{ (1+\nu)(1-2\nu)} \qquad
  \mu     = \frac{E    }{2(1+\nu)        }
\end{equation}
Furthermore $\T{I}$ is the second-order identity tensor and $\TT{I}^s$ is the fourth-order symmetrization tensor (see Appendix~\ref{sec:nomenclature}).

The connection to the deformation gradient tensor $\T{F}$ and the first Piola-Kirchhoff stress $\T{P}$ is made by the definition of the Green-Lagrange strain
\begin{equation}
  \T{E} = \tfrac{1}{2} \big( \T{F}^T \cdot \T{F} - \T{I} \big)
\end{equation}
and of the first Piola-Kirchhoff stress
\begin{equation}\label{eq:elas:S-P}
  \T{P} = \T{F} \cdot \T{S}
\end{equation}
%

\subsection{Consistent tangent}

The constitutive model is linearized in the reference configuration in accordance with Eq.~\eqref{eq:iter:varP}. Its derivation is straightforward. The first step is to linearize the stress relation in Eq.~\eqref{eq:elas:S-P}, leading to
\begin{equation}\label{eq:elas:linP}
  \delta \T{P}^T
  =
  \T{S}_{(i)} \cdot \delta \T{F}^T +
  \TT{I}^{RT} : \big( \T{F}_{(i)}  \cdot \delta \T{S} \big)
\end{equation}
wherein $\TT{I}^{RT}$ is the fourth-order right-transposed identity tensor (see Appendix~\ref{sec:nomenclature}). The constitutive model in Eq.~\eqref{eq:elas:S} is already linear. Combined with the linearization of the Green-Lagrange strain one obtains
\begin{equation}\label{eq:elas:linE}
  \delta \T{S} = \TT{C} : \big( \TT{I}^s \cdot \T{F}_{(i)}^T \big) : \delta \T{F}
\end{equation}
The expression of the tangent stiffness follows by combining (\ref{eq:elas:linP},\ref{eq:elas:linE}) to obtain
\begin{equation}\label{eq:elas:K}
  \TT{K}_{(i)}
  =
  \T{S}_{(i)} \cdot \TT{I}
  +
  \TT{I}^{RT} :
  \big( \T{F}_{(i)} \cdot \TT{C} \cdot \T{F}_{(i)}^T \big)
  : \TT{I}^{RT}
\end{equation}
whereby the right-symmetry of $\TT{C}$ has been used to absorb $\TT{I}^s$.

\subsection{Example}

\begin{figure}[htp]
  \centering
  \includegraphics[width=0.9\textwidth]{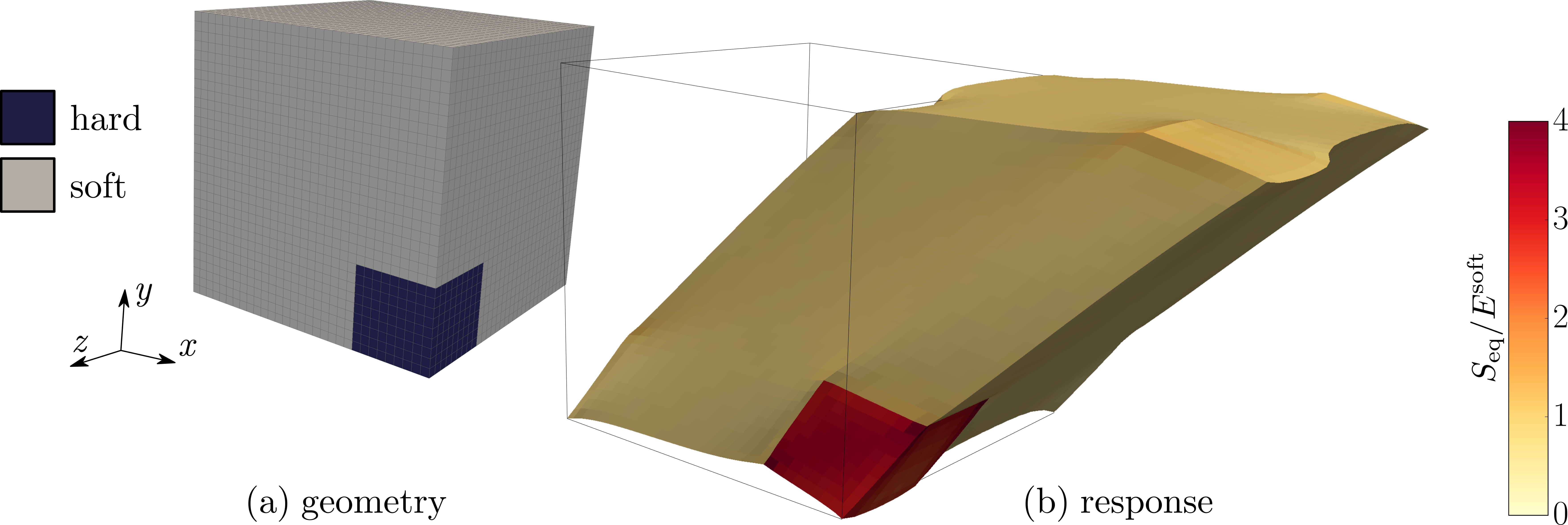}
  \caption{Simulation using the Python code of Appendix~\ref{sec:code:python}: (a) the phase distribution, and (b) the equivalent second Piola-Kirchhoff stress $S_\mathrm{eq}$ normalized with Young's modulus of the soft phase. The deformed geometry has been reconstructed using Finite Elements.}
  \label{fig:ex:lin}
\end{figure}

As an example, a three-dimensional unit-cell comprising two phases is considered -- a cubic particle with a stiffness that is $10$ times higher than that of the matrix in which it is embedded, see Figure~\ref{fig:ex:lin}(a). The unit-cell is loaded in simple shear, described by the following macroscopic deformation gradient
\begin{equation}\label{eq:elas:deform}
  \bar{\T{F}} = \T{I} + \bar{\gamma}\, \vec{e}_\mathrm{x} \vec{e}_\mathrm{y}
\end{equation}
This example has been simulated with the code given in Appendix~\ref{sec:code:python}; the parameters are $E^\mathrm{hard} = 10 E^\mathrm{soft}$ and $\nu = 0.3$. A large deformation is applied, with $\bar{\gamma} = 1$.

The example needs less than a minute to run on an ordinary desktop computer, requiring five iterations (in a single increment). The response is depicted in Figure~\ref{fig:ex:lin}(b), showing the equivalent second Piola-Kirchhoff stress $S_\mathrm{eq}$. As observed, the stress in the hard phase is significantly higher than in the soft phase, while the soft phase accommodates most of the deformation. Furthermore, even though $\T{F}$ is discontinuous across the phase boundary no artifacts of the continuous interpolation are observed at the nodes, in-spite of the phase contrast (see also \citep{Moulinec1998,Anglin2014,Zeman2016}). Note that the deformed geometry has been reconstructed for visualization purposes using Finite Elements. To this end, a mesh is used in which each element encloses one pixel. The nodal displacements of the corners are sought such that the resulting deformation gradients closely approximate the deformation gradient field $\colT{F}$ that follows from the equilibrium simulation with the FFT.

It is remarked once more that the regular grid required by the FFT can be a restriction depending on the type of problem. In the context of this example, if our goal would have been to find the homogenized response for a specific volume fraction of the inclusions, the results would be somewhat sensitive to the number of voxels used for coarse and moderately fine discretizations. They do however converge upon grid refinement.

\section{Application to Simo elasto-plasticity, defined in the current configuration}
\label{sec:simo}

\subsection{Introduction}

The elasto-plastic Simo model \citep{Simo1992a} is considered to demonstrate that a model in the current configuration can be treated in a straightforward manner, by employing the procedures developed for Finite Elements. Note that the code for this model is longer due to the more involved stress update and constitutive tangent, but its implementation is still well below 200 lines of Python code \citep{GooseFFT}.

\subsection{Constitutive model}

The model departs from the conventional multiplicative split of the deformation gradient tensor into elastic and plastic parts, $\T{F}_\mathrm{e}$ and $\T{F}_\mathrm{p}$ respectively:
\begin{equation}
  \T{F} = \T{F}_\mathrm{e} \cdot \T{F}_\mathrm{p}
\end{equation}
The Kirchhoff stress depends on the elastic part of the deformation gradient though hyperelasticity, in the form of a linear relation with the elastic logarithmic strain $\tfrac{1}{2} \ln \T{b}_\mathrm{e}$, as follows:
\begin{equation}
  \bm{\tau} = \tfrac{1}{2}\, \TT{C} : \ln \T{b}_\mathrm{e}
\end{equation}
wherein $\TT{C}$ is the elasticity tensor given in Eq.~\eqref{eq:elas:C}.

The elastic domain is bounded by the yield criterion
\begin{equation}\label{eq:simo:flow}
  \phi ( \bm{\tau} , \varepsilon_\mathrm{p} )
  =
  \tau_\mathrm{eq} - \tau_\mathrm{y} ( \varepsilon_\mathrm{p} )
  \leq 0
\end{equation}
wherein $\varepsilon_\mathrm{p}$ is the accumulated plastic strain, $\tau_\mathrm{eq}$ is the equivalent Kirchhoff stress, and $\tau_\mathrm{y}$ is the yield stress which hardens depending on $\varepsilon_\mathrm{p}$. The directions of plastic flow follow from normality. This model is implemented using an implicit time-integration scheme for the plasticity, resulting in an elastic-predictor plastic-corrector algorithm. This procedure is rather standard and can be found in several references \citep{Geers2004,Simo1992a}. Note that, by virtue of the logarithmic strain, its implementation is similar to a small-strain elasto-plastic model, available in textbooks \citep{Simo1998a,DeSouzaNeto2008}.

\subsection{Consistent tangent}

To obtain the consistent tangent we start from the linearized relation from the updated Lagrange-framework, which is obtained by linearizing about the last-known (iterative) configuration. The linearization is of the form
\begin{equation}\label{eq:simo:tangent}
  \delta \bm{\tau}
  =
  \left( -\TT{I}^{RT} \cdot \bm{\tau}_{(i)} + \TT{C}_{(i)} \right)
  : \T{F}_{(i)}^{-T} \cdot \delta \T{F}^T
\end{equation}
The first term in \eqref{eq:simo:tangent} arises from the geometric non-linearity, while the second term follows from the linearization of the constitutive model including the strain definition. For this model an explicit expression for $\TT{C}_{(i)}$ can be found in references \citep{Simo1992a,Geers2004}. Note that it involves a discrete elasto-plastic switch, which affects the convergence around the yield point.

\subsection{Pull-back to the undeformed configuration}

The pull-back of the Kirchhoff stress to the first Piola-Kirchhoff stress reads
\begin{equation}
  \T{P}_{(i)} = \bm{\tau}_{(i)} \cdot \T{F}^{-T}_{(i)}
\end{equation}
This expression is also used to pull-back the consistent tangent operator, which in the form of Eq.~\eqref{eq:simo:tangent} reads:
\begin{equation}\label{eq:pull-back}
  \delta \bm{\tau}
  =
  (\TT{K}_{\vec{x}})_{(i)}
  : \T{F}_{(i)}^{-T} \cdot \delta \T{F}^T
\end{equation}
The subscript $\vec{x}$ indicates that the tangent is defined in the last-known iterative configuration. In terms of a variation in the first Piola-Kirchhoff stress~\eqref{eq:pull-back} can be rewritten as
\begin{equation}
  \delta \T{P}^T
  =
  \big[
  \T{F}_{(i)}^{-1} \cdot
  (\TT{K}_{\vec{x}})_{(i)}
  \cdot \T{F}_{(i)}^{-T}
  \big]
  : \delta \T{F}^T
\end{equation}
Note that this pull-back does not introduce additional approximations.

\subsection{Example}

\begin{figure}[htp]
  \centering
  \includegraphics[width=1.\textwidth]{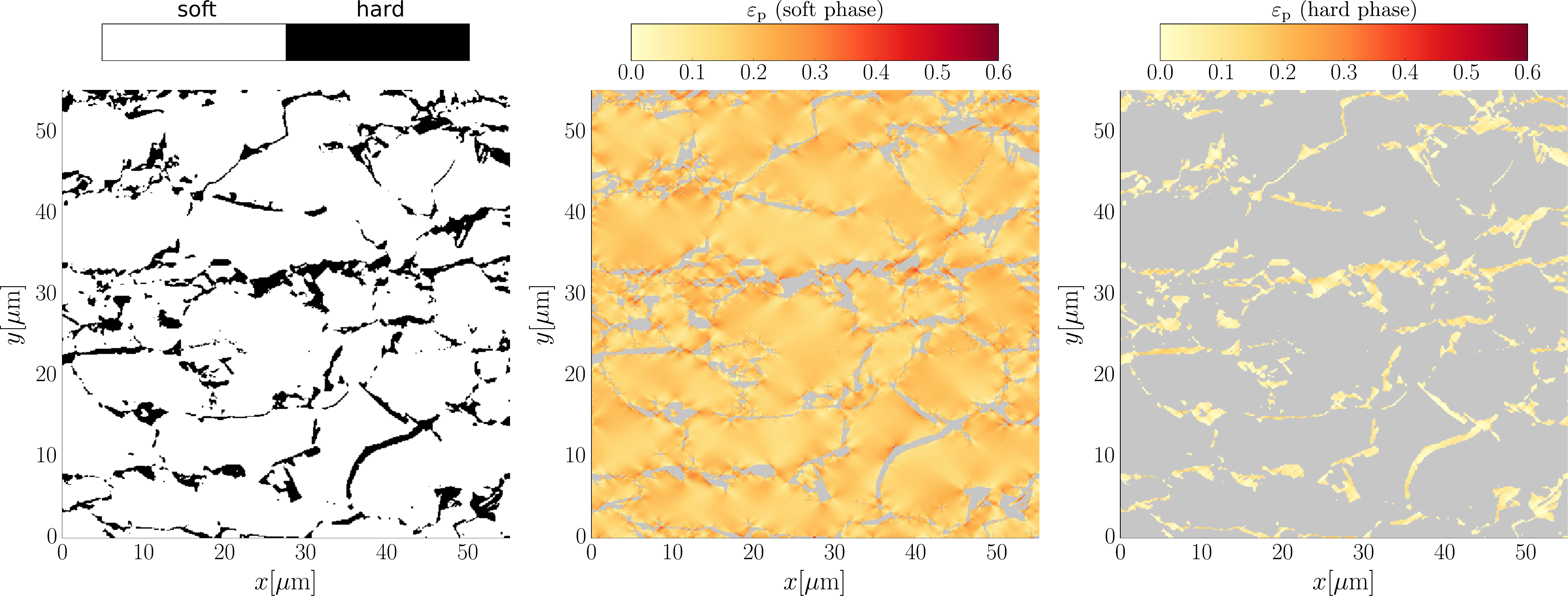}
  \caption{(a) Segmented micrograph used in the simulation. (b--c) The accumulated plastic strain $\varepsilon_\mathrm{p}$ is the soft and the hard phase respectively; for $\bar{\varepsilon} = 0.2$.}
  \label{fig:simo}
\end{figure}

To illustrate the method, a simulation is performed directly on a dual-phase steel micrograph, see Figure~\ref{fig:simo}(a), under the plane strain assumption. The micrograph has been obtained by thresholding an image of $451 \times 451$ pixels, acquired using scanning electron microscopy. For this material, the contrast in secondary electron mode is due to a small height difference between the two phases -- martensite (black in Figure~\ref{fig:simo}(a)) and ferrite (white in Figure~\ref{fig:simo}(a)) -- which results from a protocol of grinding, polishing, and etching. The two phases are assumed elastically homogeneous. Plasticity is governed by linear hardening, whereby the yield stress in \eqref{eq:simo:flow} reads
\begin{equation}
  \tau_\mathrm{y} ( \varepsilon_\mathrm{p} )
  =
  \tau_\mathrm{y0} + H \varepsilon_\mathrm{p}
\end{equation}
The parameters used for the two phases are
\begin{equation}
  \tau_\mathrm{y0}^\mathrm{soft} = \tau_\mathrm{y0}^\mathrm{hard} / \chi = 0.003 \, E
  \qquad
  H^\mathrm{soft} = H^\mathrm{hard} / \chi = 0.01 \, E
  \qquad
  \nu = 0.3
\end{equation}
where $\chi$ characterizes the phase contrast. For the example in Figure~\ref{fig:simo} we use $\chi = 2$.

Macroscopic pure shear is applied, which corresponds to the following macroscopic logarithmic strain
\begin{equation}
  \bar{\T{F}} =
  \bar{\lambda}       \vec{e}_\mathrm{x} \vec{e}_\mathrm{x} +
  1 / \bar{\lambda}\, \vec{e}_\mathrm{y} \vec{e}_\mathrm{y}
\end{equation}
where $\bar{\lambda}$ is applied in $250$ increments up to $\bar{\lambda} = 1.2$. The overall simulation takes approximately five hours to run on a desktop computer, a single increment thus takes about two minutes on average. On average $2.2$ iterations are needed for each increment: up to five around the yield point, and two in the plastic regime (in this case for $\bar{\lambda} > 1.04$).

The accumulated plastic strain response is shown in Figures~\ref{fig:simo}(b--c) for the soft and the hard phase respectively. The plastic strain is significantly higher in the soft phase than in the hard phase. The region of the localized strain is aligned with the applied shear (at $\pm 45$ degree angles) and reaches its maximum near the hard phase. With the value of $0.51$, the maximum is significantly higher than the applied macroscopic equivalent strain (of approximately $0.2$), due to the strain partitioning between the phases.

\begin{figure}[htp]
  \centering
  \includegraphics[width=1.\textwidth]{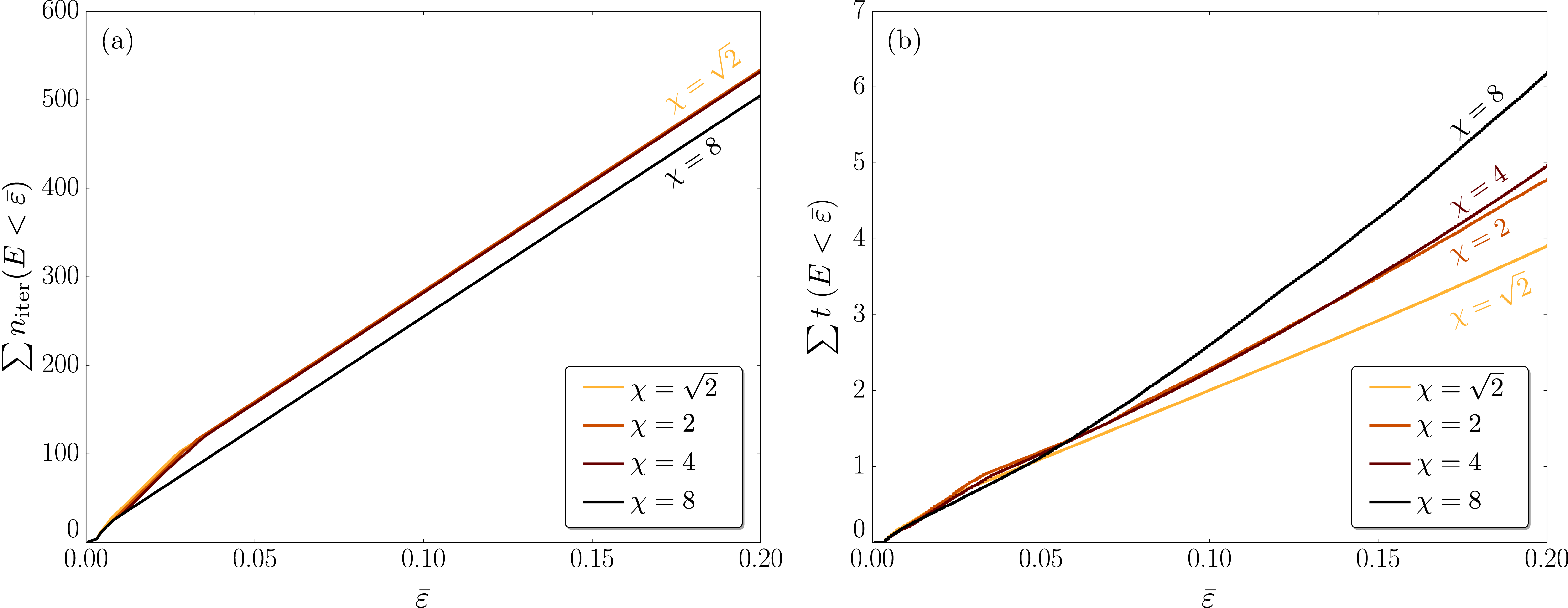}
  \caption{Numerical performance of the method for the presented example (in Figure~\ref{fig:simo}) for different phase contrasts $\chi$: (a) the cumulative number of iterations $n_\mathrm{iter}$ and (b) the cumulative runtime $t$, as a function of the applied macroscopic equivalent strain $\bar{\varepsilon}$.}
  \label{fig:simo:chi}
\end{figure}

This example finally allows us to touch upon the fact that the efficiency of the conjugate gradient solver used to solve Eq.~\eqref{eq:iter:system} depends on the conditioning of the linear system. For this example this implies that the efficiency depends on the phase contrast $\chi$. Therefore the numerical performance is compared for $\chi = \sqrt{2}, 2, 4, 8$ in Figure~\ref{fig:simo:chi}. It is observed that the Newton iterations -- characterized by the number of iterations, in Figure~\ref{fig:simo:chi}(a) -- is practically independent of $\chi$, except that convergence is faster around the yield point for $\chi = 8$. The conjugate gradient solver -- characterized by the runtime,\footnote{Note that the number of CG iterations is not directly available from the standard library routine that is used in this Python implementation. The runtime is fully equivalent here, because the number of Newton iterations is practically independent of $\chi$.} in Figure~\ref{fig:simo:chi}(b) -- however is more efficient for lower values of $\chi$, as, for that case, the entries in the constitutive $\matTT{K}$ are more homogeneous.

\section{Summary}
\label{sec:summary}

A methodology has been presented that solves large deformation mechanical equilibrium problems defined on a periodic unit-cell. It follows the standard steps used in the Finite Element Method by defining a weak form, its discretization using Galerkin's method, numerical quadrature, and finally linearization towards a Newton scheme. The differences with respect to standard Finite Elements are that (i) globally supported trigonometric polynomials are used for the interpolation and (ii) the method is formulated in terms of the nodal deformation gradient tensor and the first Piola-Kirchhoff stress, tensor, rather than nodal displacements and nodal forces. Compatibility of the solution is enforced through the projection operator in conjunction with a suitable iterative linear solver such as conjugate gradients, GMRES, Chebyshev's method, or Richardson iteration.

The following conclusions result from our developments:
\begin{itemize}
  \item The method is quite general, and is suitable for arbitrary constitutive models in a finite deformation setting.
  \item The method is memory efficient, without relying on an involved data-storage.
  \item The computational efficiency results from the use of the Fast Fourier Transform.
  \item Because of the former two facts, the method is easy to implement. This has been substantiated with a 59 lines Python code (without comments) for a 3-D finite strain problem.
  \item The constitutive level is independent of the adopted solution strategy. The same constitutive implementation can be used for a FFT-solver as well as for a conventional Finite Element solver.
\end{itemize}
%

\section*{Acknowledgments}

This research was carried out under project number M22.2.11424 in the framework of the research program of the Materials innovation institute M2i (\href{http://www.m2i.nl}{www.m2i.nl}). Jan Zeman's work was partially supported by the Czech Science Foundation under project No.~14-00420S.

\clearpage
\appendix

\section{Nomenclature}
\label{sec:nomenclature}

\newcommand{\wa}[1]{{#1}\nphantom{#1}\phantom{\T{C}   }}
\newcommand{\wb}[1]{{#1}\nphantom{#1}\phantom{C_{ijmn}}}

\def\arraystretch{1.5}
\begin{tabular}{lll}
$     \vec{a}                            $ & vector                                        & $     a_i                                       $ \\
$     \T{A}                              $ & second-order tensor                           & $     A_{ij}                                    $ \\
$     \TT{A}                             $ & fourth-order tensor                           & $     A_{ijkl}                                  $ \\
$\wa{ \T{C}   } = \T{A}^T                $ & tensor transpose                              & $\wb{ C_{ji}        } =      A_{ij}             $ \\
$\wa{ \TT{C}  } = \TT{A}^T               $ & fourth-order tensor transpose                 & $\wb{ C_{lkji}      } =      A_{ijkl}           $ \\
$\wa{ \TT{C}  } = \TT{A}^{LT}            $ & fourth-order left tensor transpose            & $\wb{ C_{jikl}      } =      A_{ijkl}           $ \\
$\wa{ \TT{C}  } = \TT{A}^{RT}            $ & fourth-order right tensor transpose           & $\wb{ C_{ijlk}      } =      A_{ijkl}           $ \\
$\wa{ \vec{c} } = \vec{a} \times  \vec{b}$ & cross product                                 & $\wb{ c_{i}         } = \epsilon_{ijk} \, a_j \, b_k$ \\
$\wa{ \T{C}   } = \vec{a} \otimes \vec{b}$ & outer product                                 & $\wb{ C_{ij}        } =      a_{i}      b_{j}   $ \\
$\wa{ \TT{C}  } = \T{A}   \otimes \T{B}  $ & outer product                                 & $\wb{ C_{ijkl}      } =      A_{ij}     B_{kl}  $ \\
$\wa{ \T{C}   } = \T{A}   \cdot   \T{B}  $ & dot product                                   & $\wb{ C_{ik}        } = \wb{ A_{ij}   } B_{jk}  $ \\
$\wa{ \TT{C}  } = \TT{A}  \cdot   \T{B}  $ & dot product                                   & $\wb{ C_{ijkm}      } = \wb{ A_{ijkl} } B_{lm}  $ \\
$\wa{ \T{C}   } = \TT{A}  :       \T{B}  $ & double dot product                            & $\wb{ C_{ij}        } = \wb{ A_{ijkl} } B_{lk}  $ \\
$\wa{ \TT{C}  } = \TT{A}  :       \TT{B} $ & double dot product                            & $\wb{ C_{ijmn}      } = \wb{ A_{ijkl} } B_{lkmn}$ \\
$\wa{ \T{C}   } = \TT{A}  \star   \T{B}  $ & convolution & $\int_{-\infty}^{\infty} \TT{A}(\vec{x}) : \T{B} (\vec{x} - \vec{y}) \;\mathrm{d}\Omega_0$ \\
$\vec{\nabla}_0 \cdot \T{A}$               & divergence (reference configuration)          & $\partial A_{ij} / \partial X_i                 $ \\
$\wa{ \T{F}   } = \big( \vec{\nabla}_0 \, \vec{x} \big)^T$ & deformation gradient tensor   & $\wb{ F_{ij}        } = \partial x_{i} / \partial X_j$ \\
$\T{I}$                                    & second order identity tensor                  & $\wb{ I_{ij}        } = \delta_{ij}             $ \\
$\TT{I}$                                   & fourth order identity tensor                  & $\wb{ I_{ijkl}      } = \delta_{il} \, \delta_{jk} $ \\
$\TT{I}^{RT}$                              & fourth order right-transposed identity tensor & $\wb{ I_{ijkl}^{RT} } = \delta_{ik} \, \delta_{jl} $ \\
$\wa{ \TT{I}^{s} } = \tfrac{1}{2} (\TT{I} + \TT{I}^{RT})$ & fourth order symmetrization tensor \\
$\colT{A}$                                 & column of tensors & $\big[ \T{A}_{1}, \T{A}_{2}, \ldots, \T{A}_{n} \big]^\mathsf{T}$ \\
$\matT{A}$                                 & matrix of tensors \\
$\colT{A}^\mathsf{T}$                      & matrix transpose  \\
\end{tabular}

\clearpage
\section{Python code}
\label{sec:code:python}

The example of Section~\ref{sec:elas} has been simulated using the Python code included in this appendix. All quantities are stored as three-dimensional matrices of second- or fourth-order tensors. This obsoletes sparse data-structures, while never storing more than needed. The implementation is rather efficient; the only somewhat costly part is the definition of the Fourier coefficients of the projection operator, $\mat{\hat{\TT{G}}}$, in lines 38--42. This is here done in this manner to enhance readability; a vectorized version is available online \citep{GooseFFT}.

What follows is the full code, interrupted by a brief explanation of the statements which follow. To run, we suggest to download the code directly from \citep{GooseFFT}.

\paragraph*{Import modules}

Load modules to extend the standard Python functionality with ``\texttt{import}''. Functions from these modules are then called using ``\texttt{modulename.''function}. This is optionally aliased, such as for ``\texttt{numpy}'' and ``\texttt{scipy.''sparse.linalg} below. NumPy is the most important module in the present context. It contains a large set of linear algebraic operations, rendering it a powerful toolkit for scientific computations.

\vspace*{1eM}%
\noindent%
\begin{minipage}{\linewidth}%
\noindent\rule{\linewidth}{0.4pt}%
\begin{lstlisting}[language=Python,firstnumber=1,firstline=1,lastline=3]
import numpy as np
import scipy.sparse.linalg as sp
import itertools
\end{lstlisting}%
\vspace*{-0.8eM}%
\rule{\linewidth}{0.4pt}%
\end{minipage}%

\paragraph*{Problem dimensions}

Define the problem dimensions: number of spatial dimensions, and number of grid points in each spatial dimension.

\vspace*{1eM}%
\noindent%
\begin{minipage}{\linewidth}%
\noindent\rule{\linewidth}{0.4pt}%
\begin{lstlisting}[language=Python,firstnumber=7,firstline=7,lastline=8]
ndim   = 3
N      = 31
\end{lstlisting}%
\vspace*{-0.8eM}%
\rule{\linewidth}{0.4pt}%
\end{minipage}%

\paragraph*{Tensor products}

Create functions to evaluate the tensor operations and products (see Appendix~\ref{sec:nomenclature}). Note that in Python, functions are defined using ``\texttt{def}'' for functions that comprise more than one line and ``\texttt{lambda}'' for single-line functions. The NumPy-function ``\texttt{einsum}'' allows the use of index notation; the comma thereby separates input arguments.

\vspace*{1eM}%
\noindent%
\begin{minipage}{\linewidth}%
\noindent\rule{\linewidth}{0.4pt}%
\begin{lstlisting}[language=Python,firstnumber=14,firstline=14,lastline=20]
trans2 = lambda A2   : np.einsum('ijxyz          ->jixyz  ',A2   )
ddot42 = lambda A4,B2: np.einsum('ijklxyz,lkxyz  ->ijxyz  ',A4,B2)
ddot44 = lambda A4,B4: np.einsum('ijklxyz,lkmnxyz->ijmnxyz',A4,B4)
dot22  = lambda A2,B2: np.einsum('ijxyz  ,jkxyz  ->ikxyz  ',A2,B2)
dot24  = lambda A2,B4: np.einsum('ijxyz  ,jkmnxyz->ikmnxyz',A2,B4)
dot42  = lambda A4,B2: np.einsum('ijklxyz,lmxyz  ->ijkmxyz',A4,B2)
dyad22 = lambda A2,B2: np.einsum('ijxyz  ,klxyz  ->ijklxyz',A2,B2)
\end{lstlisting}%
\vspace*{-0.8eM}%
\rule{\linewidth}{0.4pt}%
\end{minipage}%

\clearpage
\paragraph*{Define tensors}

Initialize grids of second- and fourth-order identity tensors are initialized to evaluate the constitutive response (see Appendix~\ref{sec:nomenclature}). Notice once more the use of the index notation in the NumPy-function ``\texttt{einsum}''.

\vspace*{1eM}%
\noindent%
\begin{minipage}{\linewidth}%
\noindent\rule{\linewidth}{0.4pt}%
\begin{lstlisting}[language=Python,firstnumber=23,firstline=23,lastline=23]
i      = np.eye(ndim)
\end{lstlisting}%
\begin{lstlisting}[language=Python,firstnumber=25,firstline=25,lastline=29]
I      = np.einsum('ij,xyz'           ,                  i   ,np.ones([N,N,N]))
I4     = np.einsum('ijkl,xyz->ijklxyz',np.einsum('il,jk',i,i),np.ones([N,N,N]))
I4rt   = np.einsum('ijkl,xyz->ijklxyz',np.einsum('ik,jl',i,i),np.ones([N,N,N]))
I4s    = (I4+I4rt)/2.
II     = dyad22(I,I)
\end{lstlisting}%
\vspace*{-0.8eM}%
\rule{\linewidth}{0.4pt}%
\end{minipage}%

\paragraph*{Projection}

Initialize the Fourier coefficients $\matTT{\hat{G}}$, cf.\ Eq.~\eqref{eq:G}. The ``\texttt{itertools}''-module allows one to loop over several indices at once. Note that in Python loops are terminated by changing the indentation level; no ``\texttt{end}''-statements are used.

\vspace*{1eM}%
\noindent%
\begin{minipage}{\linewidth}%
\noindent\rule{\linewidth}{0.4pt}%
\begin{lstlisting}[language=Python,firstnumber=34,firstline=34,lastline=36]
delta  = lambda i,j: np.float(i==j)
freq   = np.arange(-(N-1)/2.,+(N+1)/2.)
Ghat4  = np.zeros([ndim,ndim,ndim,ndim,N,N,N])
\end{lstlisting}%
\begin{lstlisting}[language=Python,firstnumber=38,firstline=38,lastline=42]
for i,j,l,m in itertools.product(range(ndim),repeat=4):
    for x,y,z    in itertools.product(range(N),   repeat=3):
        q = np.array([freq[x], freq[y], freq[z]])
        if not q.dot(q) == 0:
            Ghat4[i,j,l,m,x,y,z] = delta(i,m)*q[j]*q[l]/(q.dot(q))
\end{lstlisting}%
\vspace*{-0.8eM}%
\rule{\linewidth}{0.4pt}%
\end{minipage}%

\vspace*{1eM}\noindent
Create functions to evaluate the operations $\matTT{G} : \bullet$ and $\matTT{G} : \matTT{K}_{(i)} : \bullet$ using a sequence of the Fourier transform $\mathcal{F}$, the double tensor contraction with the Fourier coefficients $\mat{\hat{\TT{G}}}$, and the inverse Fourier transform $\mathcal{F}^{-1}$. See Equations~(\ref{eq:weak:discr:FFT-intermediate},\ref{eq:weak:discr:FFT}). Note that ``\texttt{fftshift}'', and its inverse, have been used to align the zero frequency with the center of the grid.

\vspace*{1eM}%
\noindent%
\begin{minipage}{\linewidth}%
\noindent\rule{\linewidth}{0.4pt}%
\begin{lstlisting}[language=Python,firstnumber=45,firstline=45,lastline=46]
fft    = lambda x  : np.fft.fftshift(np.fft.fftn (np.fft.ifftshift(x),[N,N,N]))
ifft   = lambda x  : np.fft.fftshift(np.fft.ifftn(np.fft.ifftshift(x),[N,N,N]))
\end{lstlisting}%
\begin{lstlisting}[language=Python,firstnumber=49,firstline=49,lastline=51]
G      = lambda A2 : np.real( ifft( ddot42(Ghat4,fft(A2)) ) ).reshape(-1)
K_dF   = lambda dFm: trans2(ddot42(K4,trans2(dFm.reshape(ndim,ndim,N,N,N))))
G_K_dF = lambda dFm: G(K_dF(dFm))
\end{lstlisting}%
\vspace*{-0.8eM}%
\rule{\linewidth}{0.4pt}%
\end{minipage}%

\clearpage
\paragraph*{Problem definition and constitutive model}

Create a grid of scalars for the bulk modulus $K$ and the shear modulus $G$. For the example of Section~\ref{sec:elas} a cubic particle is embedded in an otherwise homogeneous matrix, see also Figure~\ref{fig:ex:lin}.

\vspace*{1eM}%
\noindent%
\begin{minipage}{\linewidth}%
\noindent\rule{\linewidth}{0.4pt}%
\begin{lstlisting}[language=Python,firstnumber=56,firstline=56,lastline=56]
phase  = np.zeros([N,N,N]); phase[-9:,:9,-9:] = 1.
\end{lstlisting}%
\begin{lstlisting}[language=Python,firstnumber=58,firstline=58,lastline=60]
param  = lambda M0,M1: M0*np.ones([N,N,N])*(1.-phase)+M1*np.ones([N,N,N])*phase
K      = param(0.833,8.33)
mu     = param(0.386,3.86)
\end{lstlisting}%
\vspace*{-0.8eM}%
\rule{\linewidth}{0.4pt}%
\end{minipage}%

\vspace*{1eM}\noindent
Define a function to evaluate the constitutive model (for the entire grid). I.e.\ calculate the first Piola-Kirchhoff stress tensors $\colT{P}_{(i)}$ and tangent stiffness tensors $\matTT{K}_{(i)}$ based on the deformation gradients $\colT{F}_{(i)}$, see Equations~(\ref{eq:elas:S}--\ref{eq:elas:K}).

\vspace*{1eM}%
\noindent%
\begin{minipage}{\linewidth}%
\noindent\rule{\linewidth}{0.4pt}%
\begin{lstlisting}[language=Python,firstnumber=63,firstline=63,lastline=68]
def constitutive(F):
    C4 = K*II+2.*mu*(I4s-1./3.*II)
    S  = ddot42(C4,.5*(dot22(trans2(F),F)-I))
    P  = dot22(F,S)
    K4 = dot24(S,I4)+ddot44(ddot44(I4rt,dot42(dot24(F,C4),trans2(F))),I4rt)
    return P,K4
\end{lstlisting}%
\vspace*{-0.8eM}%
\rule{\linewidth}{0.4pt}%
\end{minipage}%

\paragraph*{Newton iterations}

Initialize a stress-free state with $\colT{P} = \colT{O}$ and $\colT{F} = \colT{I}$, both as a grid of second order tensors.

\vspace*{1eM}%
\noindent%
\begin{minipage}{\linewidth}%
\noindent\rule{\linewidth}{0.4pt}%
\begin{lstlisting}[language=Python,firstnumber=73,firstline=73,lastline=74]
F     = np.array(I,copy=True)
P,K4  = constitutive(F)
\end{lstlisting}%
\vspace*{-0.8eM}%
\rule{\linewidth}{0.4pt}%
\end{minipage}%

\vspace*{1eM}\noindent
Apply macroscopic deformation gradient of Equation~\eqref{eq:elas:deform}, and define the right-hand side of the linear system to distribute $\Delta \colT{\bar{F}}$ over the grid using Equation~\eqref{eq:bc:lin}. Note that the order of operations is slightly different than in Algorithm~\ref{fig:algorithm} to minimize code duplication, as is more natural for the implementation.

\vspace*{1eM}%
\noindent%
\begin{minipage}{\linewidth}%
\noindent\rule{\linewidth}{0.4pt}%
\begin{lstlisting}[language=Python,firstnumber=77,firstline=77,lastline=77]
DbarF = np.zeros([ndim,ndim,N,N,N]); DbarF[0,1] += 1.0
\end{lstlisting}%
\begin{lstlisting}[language=Python,firstnumber=80,firstline=80,lastline=83]
b     = -G_K_dF(DbarF)
F    +=         DbarF
Fn    = np.linalg.norm(F)
iiter = 0
\end{lstlisting}%
\vspace*{-0.8eM}%
\rule{\linewidth}{0.4pt}%
\end{minipage}%

\vspace*{1eM}\noindent
Iterate towards equilibrium using the standard Newton algorithm. To this end the constitutive response is evaluated to calculate the residual and the tangent ($\colT{P}_{(i)}$ and $\matTT{K}_{(i)}$). A linear system is solved using conjugate gradients to obtain the iterative update $\delta \colT{F}$, which is then added to the nodal degrees-of-freedom. See Equations~(\ref{eq:iter:F}--\ref{eq:iter:system}), and Algorithm~\ref{fig:algorithm}.

\vspace*{1eM}%
\noindent%
\begin{minipage}{\linewidth}%
\noindent\rule{\linewidth}{0.4pt}%
\begin{lstlisting}[language=Python,firstnumber=86,firstline=86,lastline=96]
while True:
    dFm,_ = sp.cg(tol=1.e-8,
      A = sp.LinearOperator(shape=(F.size,F.size),matvec=G_K_dF,dtype='float'),
      b = b,
    )
    F    += dFm.reshape(ndim,ndim,N,N,N)
    P,K4  = constitutive(F)
    b     = -G(P)
    print('%10.2e'%(np.linalg.norm(dFm)/Fn))
    if np.linalg.norm(dFm)/Fn<1.e-5 and iiter>0: break
    iiter += 1
\end{lstlisting}%
\vspace*{-0.8eM}%
\rule{\linewidth}{0.4pt}%
\end{minipage}%

%
\clearpage
\section{Even-sized grid}
\label{sec:non-odd-grid}

When the problem is approximated with trigonometric polynomials using an odd-sized grid, in the final solution the deformation gradient is compatible and the stress is equilibrated. For even-sized grids, both conditions cannot be satisfied at the same time, which is caused by the Nyquist frequencies (see \cite{Vondrejc2015a} for a detailed elaboration in a scalar elliptic setting). When the values of the projection are set to zero for the Nyquist frequencies, Eq.~\eqref{eq:G} is replaced by
\begin{equation}
(\hat{G}_{kk})_{ijlm} (\vec{q}_k) =
\begin{cases}
  \displaystyle
  0
  & \displaystyle\vphantom{\frac{\vec{\xi}}{\vec{\xi}}}
  \text{for }\vec{q}_k = \vec{0} \text{ and when }\vec{q}_k\text{ has Nyquist frequency}
  \\
  \displaystyle
  \frac{
    \delta_{im}\, \xi_j(\vec{q}_k)\, \xi_l(\vec{q}_k)
  }{
    \| \vec{\xi} \|^2
  }
  & \text{otherwise}
\end{cases}
\end{equation}
to recover a compatible deformation gradient. When we set it as identity for Nyquist frequencies, the equilibrated stress can be obtained. For the example, the former approach is implemented in Python as follows:

\vspace*{1eM}%
\noindent%
\begin{minipage}{\linewidth}%
\noindent\rule{\linewidth}{0.4pt}%
\begin{lstlisting}[language=Python,firstnumber=34,firstline=34,lastline=36]
delta  = lambda i,j: np.float(i==j)
freq   = np.arange(-N/2.,+N/2.)
Ghat4  = np.zeros([ndim,ndim,ndim,ndim,N,N,N])
\end{lstlisting}%
\begin{lstlisting}[language=Python,firstnumber=38,firstline=38,lastline=42]
for i,j,l,m in itertools.product(range(ndim),repeat=4):
    for x,y,z    in itertools.product(range(1,N), repeat=3):
        q = np.array([freq[x], freq[y], freq[z]])
        if not q.dot(q) == 0:
            Ghat4[i,j,l,m,x,y,z] = delta(i,m)*q[j]*q[l]/(q.dot(q))
\end{lstlisting}%
\vspace*{-0.8eM}%
\rule{\linewidth}{0.4pt}%
\end{minipage}%
\vspace*{1eM}%

In practice the values of the projection are set to zero for the Nyquist frequencies, since it is generally not a concern to satisfy equilibrium in an approximate manner only. In fact, in confronting the simulation with the analytical solution by \citet{Kabel2014} a perfect agreement was found. Also for the example of Section~\ref{sec:simo} the results were identical using an even grid (with one pixel more to each side).

\section{Projection operator}
\label{sec:projection}

In order to explain the rationale behind the construction of the projection operator in Eq.~\eqref{eq:G}, we first observe that, in the Fourier space, it admits the expression
\begin{align}\label{eq:G_red}
\hat{G}_{ijlm}(\vec{q}) = \delta_{im} \, \hat{g}_{jl} (\vec{q})
\quad \text{ with } \quad
\hat{\T{g}}(\vec{q})
=
\begin{cases}
  \displaystyle
  \T{0}
  &\text{for }\vec{q} = \vec{0}
  \\
  \displaystyle
  \frac{\vec{\xi}(\vec{q})}{\| \vec{\xi}(\vec{q}) \|}
  \otimes
  \frac{\vec{\xi}(\vec{q})}{\| \vec{\xi}(\vec{q}) \|}
  & \text{otherwise}
\end{cases}
\end{align}
where the scaled frequencies read $\xi_i = q_i / L_i$. Now, utilizing the standard calculus associated with the Fourier transform (see e.g.~\cite{Rudin1987}), we rewrite the convolution $\T{A} = \TT{G} \star \T{B}$ as follows:
\begin{align}
\hat{A}_{ij}(\vec{q})
=
\hat{G}_{ijlm}(\vec{q}) \,
\hat{B}_{ml}(\vec{q})
=
\delta_{im} \,
\hat{g}_{jl}(\vec{q}) \,
\hat{B}_{ml}(\vec{q})
\end{align}
This reveals that the convolution of $\TT{G}$ with $\T{B}$  in fact represents convolution of $\T{g}$ with all rows of $\T{B}$. Therefore, it suffices to concentrate on the relationship between an arbitrary row of $\T{A}$, say $\vec{a}$, and the corresponding row of $\T{B}$, say $\vec{b}$:
\begin{equation}
\vec{a} = \T{g} \star \vec{b}
\end{equation}

To show that $\vec{a}$ can be obtained as the gradient of a scalar potential, we need to verify that its rotation~(or curl) vanishes, i.e.~$\vec{\nabla}_0 \times \vec{a} = \vec{0}$. In the Fourier domain, the curl $\vec{\nabla}_0 \times \vec{a}$ transforms into
\begin{align}
\vec{\xi}(\vec{q}) \times \hat{\vec{a}}( \vec{q} )
=
\vec{\xi}(\vec{q}) \times \left( \hat{\T{g}}(\vec{q}) \cdot \hat{\vec{b}}(\vec{q}) \right)
=
\begin{cases}
  \displaystyle
  \vec{0}
  & \text{for }\vec{q} = \vec{0}
  \\
  \displaystyle
  \vec{\xi}(\vec{q}) \times
  \frac{\vec{\xi}(\vec{q})}{\| \vec{\xi}(\vec{q}) \|}
  \frac{\vec{\xi}(\vec{q}) \cdot \hat{\vec{b}}(\vec{q})}{\| \vec{\xi}(\vec{q}) \|}
  & \text{otherwise}
\end{cases}
\end{align}
Because $\hat{\T{g}}$ projects the Fourier coefficient $\hat{\vec{b}}(\vec{q})$ to the direction of the scaled frequency $\vec{\xi}(\vec{q})$, also the second term vanishes (for $\vec{q} \neq \vec{0}$).

We have shown that all projected fields are compatible, but we still need to demonstrate that the range of $\T{g}$ coincides with the \emph{whole} set of compatible fields, not only with a subset. To this purpose, we express a compatible field $\vec{b}$ in terms of a scalar potential $f$ via $\vec{b} = \vec{\nabla}_0 f$, which in Fourier space reads $\hat{\vec{b}}(\vec{q}) = \vec{\xi}(\vec{q}) \hat{f}(\vec{q})$. We thus find that
\begin{align}
\hat{\T{g}}(\vec{q}) \cdot \hat{\vec{b}}(\vec{q})
=
\begin{cases}
  \displaystyle
  \vec{0}
  & \text{for }\vec{q} = \vec{0}
  \\
  \displaystyle
  \frac{\vec{\xi}(\vec{q})}{\| \vec{\xi}(\vec{q}) \|}
  \frac{\vec{\xi}(\vec{q}) \cdot \vec{\xi}(\vec{q})}{\| \vec{\xi}(\vec{q}) \|}
  \hat{f}(\vec{q})
  & \text{otherwise}
\end{cases}
\end{align}
or simply
\begin{equation}
\hat{\T{g}}(\vec{q}) \cdot \hat{\vec{b}}(\vec{q})
=
\vec{\xi}(\vec{q}) \hat{f}(\vec{q})
=
\hat{\vec{b}}(\vec{q})
\end{equation}
Altogether, this shows that operator $\TT{G}$ maps square-integrable second-order tensorial fields onto the zero-mean compatible ones.

\bibliography{library}

\end{document}